# Group invariant inferred distributions via noncommutative probability

B. Heller[1] and M. Wang[2],*

*Illinois Institute of Technology and University of Chicago*

**Abstract:** One may consider three types of statistical inference: Bayesian, frequentist, and group invariance-based. The focus here is on the last method. We consider the Poisson and binomial distributions in detail to illustrate a group invariance method for constructing inferred distributions on parameter spaces given observed results. These inferred distributions are obtained without using Bayes' method and in particular without using a joint distribution of random variable and parameter. In the Poisson and binomial cases, the final formulas for inferred distributions coincide with the formulas for Bayes posteriors with uniform priors.

## 1. Introduction

The purpose of this paper is to construct a probability distribution on the parameter space given an observed result (*inferred* distribution) in the case of a discrete random variable with a continuous parameter space using group theoretic methods. We present two examples, the Poisson and binomial distributions. From the point of view of posterior Bayes distributions, these group theoretic methods lead to uniform non-informative prior distributions on canonical parameters for both the Poisson and binomial cases. Alternatively, posterior distributions are obtained here from group theory alone without explicitly using Bayes method.

The construction of inferred probability distributions by non-Bayesian methods has a long history beginning with Fisher's fiducial method of inference. The use of group theoretic methods to construct pivotal functions also has a long history as introduced in Fraser (1961) and amplified by many others since then.

Briefly, the group theoretic or "invariance" method of inference has operated from a context in which a group $G$ acts upon both the parameter space and the sample space. Consider the description as given in Eaton (1989). Let $(\mathbf{X}, \mathcal{B})$ represent a given measurable space associated with random variable $X$ having probability distribution $Q_0$. Assume that there is an action of group $G$ on the sample space $\mathbf{X}$ of the random variable. There is an induced action of $G$ on probability distributions. Thus, if $X$ has probability distribution $Q_0$ then define $gQ_0$ as the probability distribution of the random variable $gX$. (This is done similarly in Fraser (1961) for probability distributions of sufficient statistics.) Then consider the collection of probability distributions $\{gQ_0 | g \in G\}$. If $G$ is a Lie group (i.e. parameterized) then we have a collection of probability distributions indexed by the group parameters and an associated invariant measure on the parameter space.

*Corresponding author.
[1]Department of Applied Mathematics, Illinois Institute of Technology, Chicago, IL 60616, e-mail: heller@iit.edu
[2]Department of Statistics, University of Chicago, Chicago, IL 60637, e-mail: meiwang@galton.uchicago.edu

The salient feature of this (essentially) pivotal method is an isomorphism between three entities: the group $G$, the parameter space, and the statistic sample space $\mathbf{X}$. Clearly it is not applicable in the case of discrete distributions with continuous parameter spaces.

Group invariance methods have also been used to obtain reference priors for Bayesian posterior distributions. A comprehensive review on the selection of prior distributions is given in Kass and Wasserman (1996). In their section on invariance methods, description is given in which the group $G$ acts on both the parameter space and the sample space as outlined above. Again, those methods are not applicable to discrete distributions with continuous parameter spaces.

The unique contribution made in this paper is a group theoretic invariant method of inference which is indeed applicable in the discrete case. Also, we describe two examples, namely, the Poisson and binomial distributions. In this method, the (chosen) group $G$ acts on the parameter space but not necessarily on the sample space and we do not construct a pivotal function. Yet the group theory still leads us to an inferred distribution on the parameter space given an observed value of the random variable. Also, from the Bayesian point of view, our group invariance method provides reference priors in the case of discrete random variables. The key to the method is that a group is used to construct the requisite family in the first place. Then the group theory allows us to reverse directions to construct the inferred distribution on the parameter space. Technically this inference is possible due to the generalized spectral theorem.

The technical constructions of probability distributions given in this paper stem from some methods used in quantum physics which are used for purposes other than those described here. We use technical approaches related to four types of concepts. One relates to the idea basic to quantum physics of "non-commutative" probability as described in Parthasarathy (1992) and Whittle (1992). A second concept, so-called "covariant probability operator-valued measures" is used in what may be described as statistical design problems in communication theory such as those found in Holevo (2001), Helstrom (1976), and Busch, Grabowski and Lahti (1995). The third concept, "coherent states", is described in Perelomov (1986) from a strictly group theoretic point of view and more generally in Ali, Antoine and Gazeau (2000). The fourth type of material is group representation theory itself as given in Vilenkin (1968).

It should be noted that some statisticians are becoming interested in quantum physics from the point of view of how one should deal with quantum data. An overview of quantum theory and the relationship to statistical methods for dealing with quantum data is given in Malley and Hornstein (1993) and in Barndorff-Nielsen, Gill and Jupp (2003). Explanations of quantum theory and its relationship to statistical problems are outlined in the works of Helland, for example Helland (1999, 2003a, 2003b). However, in this paper we are not dealing with problems of quantum data analysis. We are simply using some technical methods which appear in the quantum physics literature as well as in the group theory literature for our own purposes.

### 1.1. Noncommutative probability distributions

The probability distributions we consider are obtained in a different manner than those of classical probability. In Parthasarathy (1992), the difference is explained in the following manner. Suppose that we consider the expectation $\mathbb{E}[Y]$ of a real valued discrete random variable. For example, suppose possible value $y_i$ has probability $p_i$, for $i = 1.2, \ldots, k$. One can express the expectation in terms of the trace of the product of two $k \times k$ diagonal matrices, $S$ and $O$:

$$\mathbb{E}[Y] = \sum_{i=1}^{k} p_i y_i = trace \left[ \begin{pmatrix} p_1 & 0 & 0 & \cdots & 0 \\ 0 & p_2 & 0 & \cdots & 0 \\ \vdots & \vdots & \vdots & \cdots & \vdots \\ 0 & 0 & 0 & \cdots & p_k \end{pmatrix} \begin{pmatrix} y_1 & 0 & 0 & \cdots & 0 \\ 0 & y_2 & 0 & \cdots & 0 \\ \vdots & \vdots & \vdots & \cdots & \vdots \\ 0 & 0 & 0 & \cdots & y_k \end{pmatrix} \right] = tr(SO).$$

In this case since the two matrices are diagonal, they are commutative. However, noncommutative matrices (or more generally, linear operators) may be used to construct expectations.

We begin by showing how to construct noncommutative probability distributions. From there we go on to generate *families* of probability distributions, and finally, we construct *inferred distributions* for the Poisson and binomial families.

We conceive of a "random experiment" as having two parts. The "input" status is represented by a linear, bounded, Hermitian, positive, trace-one operator $S$ called a *state* operator. For example, if one were tossing a coin, the bias of the coin would be represented by a state operator; loosely speaking, the *state* of the coin. The measurement process (discrete or continuous) or "output" is represented by a linear self-adjoint operator, $O$, called an *observable* or outcome operator. So that, if one tossed the coin ten times, the measurement process would be to count the number of heads. These linear operators act in a complex separable Hilbert space $\mathcal{H}$ with inner product $(\cdot, \cdot)$, which is linear in the second entry and complex conjugate linear in the first entry.

Since the observable operator is self-adjoint, it has a real spectrum. We shall consider cases where the spectrum is either all discrete or all continuous. Although operators in a Hilbert space seem far removed from a probability distribution over possible results of an experiment, the relationship is made in the following manner:

(i) The set of possible (real) results of measurement is the spectrum of the observable operator $O$. (So, in the coin tossing experiment, $O$ would have a discrete spectrum: $\{0, 1, 2, 3, 4, 5, 6, 7, 8, 9, 10\}$.)
(ii) The expected value for those results, using state operator $S$, is given by $trace(SO)$. See Whittle (1992) and Helland (2003a, b).

In order to obtain a probability distribution, the theory then depends upon the spectral theorem for self-adjoint operators. To each self-adjoint $O$, is associated a unique set of projection operators $\{\mathcal{E}(B)\}$, for any real Borel set $B$ such that

$$\mathcal{P}\{result \in B \text{ when the state operator is } S\} = trace(S \, \mathcal{E}(B)).$$

This set of of projection operators is called the *spectral measure* or the *projection-valued* (PV) measure associated with the self-adjoint operator $O$. A rigorous definition of PV measure is given in Section 2.5.

There are certain kinds of state operators that are simple to manipulate. They are the projectors onto one-dimensional subspaces spanned by unit vectors $\varphi$ in the Hilbert space $\mathcal{H}$. Since each such projection operator is identified by a unit vector in $\mathcal{H}$, the unit vector itself is called a *vector state*. In this case, the trace formula becomes simplified to an inner product: $trace(S\mathcal{E}(B)) = (\varphi, \mathcal{E}(B)\varphi)$, where $S$ is the projector onto the one-dimensional subspace spanned by unit vector $\varphi$.

Note that if unit vector $\varphi$ is multiplied by a scalar $\epsilon$ of unit modulus, we obtain the same probability distribution as with the vector $\varphi$ itself. Thus we distinguish between a single unit vector $\varphi$ and the equivalence class of unit vectors $\{\epsilon\varphi\}$ of which $\varphi$ is a representative. We use the words *vector state* or just *state* to refer to an *arbitrary* representative of a unit vector equivalence class. Thus since, for complex number $\epsilon$ of unit modulus,

$$\mathcal{P}\{result \in B \text{ when the state is } \varphi\} = \mathcal{P}\{result \in B \text{ when the state is } \epsilon\varphi\},$$

we take $\varphi$ and $\epsilon\varphi$ to be the same *state* even though they are not the same vectors.

From now on we reserve the use of the word *state* for vector states as described above. To designate a state which is an operator, as opposed to a vector, we use the phrase *state operator*.

### 1.2. Discrete probability distributions

Consider the case where the spectrum of $O$ is purely discrete and finite, consisting of eigenvalues $\{y_i\}$. Then the eigenvectors $\{\eta_i\}$ of $O$ form a complete orthonormal basis for the Hilbert space $\mathcal{H}$. In the infinite case, the Hilbert space is realized as an $\ell^2$ space of square summable sequences. When the state is $\varphi$, the probability of obtaining result $y_i$ is given by $(\varphi, \mathcal{E}(\{y_i\})\varphi)$, where $\mathcal{E}(\{y_i\})$ is the projection onto the subspace spanned by the eigenvectors of the eigenvalue $y_i$.

In particular, when the spectrum is *simple*, that is, when there is only one eigenvector $\eta_i$ for each eigenvalue $y_i$,

(1.2.1) $$\mathcal{P}\{result = y_i \text{ when the state is } \varphi\} = |(\varphi, \eta_i)|^2.$$

In order to present examples, we must first decide where to start. The natural method in the performance of statistical inference is to start with a statistical model (say a parametric family of probability distributions) pertaining to the particular physical properties of a given random experiment. Then, perhaps, one may construct posterior distributions on the parameter space based upon observed results. However, here we attempt to construct prototype families for which the inference procedures that we illustrate below, can be put in place.

Thus instead of starting with a statistical model for a particular situation, we start with an observable self-adjoint operator. As this paper progresses, it will become clear that selections of observables and families of states stem primarily from the selection of a Lie algebra. In this section, however, we consider an example of a PV measure by starting with a given observable operator in the case where its spectrum is discrete and, in fact, finite.

**Example 1.2.** Consider an experiment with three possible results $1, 0, -1$. Suppose the observable operator is

$$O = \begin{pmatrix} 1 & 0 & 0 \\ 0 & 0 & 0 \\ 0 & 0 & -1 \end{pmatrix}.$$

Note that $O$ is Hermitian. The eigenvalues of $A$ are $1, 0, -1$, and corresponding eigenvectors are

$$\eta_1 = \begin{pmatrix} 1 \\ 0 \\ 0 \end{pmatrix}, \quad \eta_2 = \begin{pmatrix} 0 \\ 1 \\ 0 \end{pmatrix}, \quad \eta_3 = \begin{pmatrix} 0 \\ 0 \\ 1 \end{pmatrix}.$$

Once the measurement is represented by a self-adjoint operator $O$ whose eigenvectors serve as a basis for the Hilbert space, then the probability distribution is determined by the choice of state.

**Part(a).** Consider the unit vector $\xi = \frac{1}{\sqrt{14}} \begin{pmatrix} 1 \\ 2 \\ 3i \end{pmatrix}$. Using (1.2.1) we have,

$$\mathcal{P}_\xi(result = 1) = |(\eta_1, \xi)|^2 = \frac{1}{14},$$
$$\mathcal{P}_\xi(result = 0) = |(\eta_2, \xi)|^2 = \frac{4}{14},$$
$$\mathcal{P}_\xi(result = -1) = |(\eta_3, \xi)|^2 = \frac{9}{14}.$$

**Part(b).** Consider the unit vector $\psi_0 = \frac{1}{2} \begin{pmatrix} -i \\ \sqrt{2} \\ i \end{pmatrix}$.

The probabilities for results $1, 0, -1$ are $\frac{1}{4}, \frac{2}{4}, \frac{1}{4}$ respectively. We see here how the choice of the state determines the probability distribution.

**Part(c).** Suppose the experiment has rotational symmetry and the probabilistic model does not change under rotation in three-dimensional space. Consider a family of states corresponding to points on the unit sphere indexed by angles $\beta$ and $\theta$ where $0 \leq \beta < 2\pi, 0 \leq \theta < \pi$. Let

$$\psi_{\beta,\theta} = \begin{pmatrix} e^{-i\beta} \cos^2 \frac{\theta}{2}, & \frac{1}{\sqrt{2}} \sin\theta, & e^{i\beta} \sin^2 \frac{\theta}{2} \end{pmatrix}^T.$$

Then $\mathcal{P}_{\psi_{\beta,\theta}}(result = 1) = |(\eta_1, \psi_{\beta,\theta})|^2 = \cos^4 \frac{\theta}{2},$

$$\mathcal{P}_{\psi_{\beta,\theta}}(result = 0) = |(\eta_2, \psi_{\beta,\theta})|^2 = \frac{1}{2}\sin^2\theta = 2\sin^2\frac{\theta}{2}\cos^2\frac{\theta}{2},$$
$$\mathcal{P}_{\psi_{\beta,\theta}}(result = -1) = |(\eta_3, \psi_{\beta,\theta})|^2 = \sin^4 \frac{\theta}{2}.$$

Relabel the possible values: $1, 0, -1$ as $0, 1, 2$, and let $p = \sin^2 \frac{\theta}{2}$. Then this family becomes the binomial distribution with $n = 2$.

### 1.3. Continuous probability distributions

In the case where the observable self-adjoint operator $O$ has a purely continuous *simple* spectrum $sp(O)$, (that is, there exists a vector $\psi_0$ in the domain of $O$ such that finite linear combinations of $O^n \psi_0$ are dense in the domain), the Hilbert space is realized as an $\mathcal{L}^2(sp(O), \mu)$ space of complex-valued square integrable functions of a real variable $x$ with inner product

$$(\psi(x), \phi(x)) = \int_{sp(O)} \psi(x)^* \phi(x) \mu(dx),$$

for some finite measure $\mu$ with support $sp(O)$, where * indicates complex conjugate. From the spectral theorem (Beltrametti and Cassinelli (1981)), we have the result that self-adjoint operator $O$ determines a unique projection-valued (PV) measure $\{\mathcal{E}(B)\}$ for real Borel sets $B$. In that case, integrating with respect to the PV measure $\mathcal{E}(B)$, we have formal operator equations:

(i) $\int_{sp(O)} \mathcal{E}(dx) = I.$

This should be understood in the sense that $(\psi(x), \mathcal{E}(dx)\psi(x))$ is a probability measure on $sp(O), \forall \psi \in$ the domain of $O$, since $(\psi(x), \mathcal{E}(B)\psi(x))$ is well defined on any Borel set $B \subseteq sp(O)$.

(ii) $O = \int_{sp(O)} x \mathcal{E}(dx).$

It follows that for certain operator functions $f(O)$, we have

$$f(O) = \int_{sp(O)} f(x) \mathcal{E}(dx).$$

In particular, let $\chi_B(x)$ be the characteristic function for Borel set $B$ and let the corresponding operator function be designated as $\chi_B(O)$. Then

$$\chi_B(O) = \int_{sp(O)} \chi_B(x) \mathcal{E}(dx).$$

For vector $\xi$ in the domain of $O$,

$$(\xi, \chi_B(O)\xi) = \int_{sp(O)} \chi_B(x) \left(\xi(x), \mathcal{E}(dx)\xi(x)\right).$$

When the Hilbert space is constructed in this manner, we say the particular Hilbert space realization is "diagonal in $O$" or is "$O$-space". In that case, the $O$ operator is the "multiplication" operator thus: $O\xi(x) = x\xi(x)$, (which explains why the spectral measure for Borel set $B$ is just the characteristic function of that set).

In the diagonal representation of the Hilbert space, since the projection operators $\{\mathcal{E}(B)\}$ are simply the characteristic function $\chi_B(O)$ for that Borel set, we have a simplified form for the probability formula. For unit vector $\psi$ in the domain of $O$:

$$\mathcal{P}_\psi(O \text{ result } \in B) = (\psi, \chi_B(O)\psi) = \int_B |\psi(x)|^2 \mu(dx).$$

Note that, in this $O$-diagonal space, the probability distribution is determined by the choice of state $\psi$.

It is possible to have spectral measures associated with operators which are not self-adjoint. Normal operators also uniquely determine spectral measures but the spectrum might be complex. Subnormal operators are associated with spectral measures in which the operators $\mathcal{F}(B)$ for complex Borel set $B$, are not projection operators but are positive operators. We will be using spectral measures of this sort, called "positive-operator-valued" (POV) measures (ref. Section 2.5), instead of projection-valued (PV) measures when we consider probability distributions on parameter spaces.

**Example 1.3.** We consider the self-adjoint operator $Q$ where $Q\psi(x) = x\psi(x), \psi \in \mathcal{L}^2(sp(Q), dx) \equiv \mathcal{H}$. The Hilbert space $\mathcal{H}$ is diagonal in $Q$, which represents the measurement of one-dimensional position in an experiment. The spectrum $sp(Q)$ of $Q$ is the whole real line $\mathbb{R}$.

We choose a state (function of $x$), $\psi(x) = \frac{e^{-x^2/4\sigma^2}}{(2\pi\sigma^2)^{1/4}}$ for $\sigma > 0$. Then

$$\mathcal{P}_\psi\{position \in B\} = (\psi(x), \mathcal{E}(B)\psi(x)) = (\psi(x), \chi_B(Q)\psi(x)) = \int_B |\psi(x)|^2 dx.$$

Thus, the probability density function for the distribution is the modulus squared of $\psi(x)$ which is the normal density function with zero mean and variance $= \sigma^2$.

### *1.4. Groups and group representations*

We will consider groups, transformation groups, Lie groups, Lie algebras, and representations of them by linear operators in various linear spaces.

#### *1.4.1. Group representations*

Let $G$ be a group with elements denoted by $g$. Let $T$ denote a linear operator on a linear space $\mathcal{H}$, a complex Hilbert space. If, to every $g \in G$, there is assigned a linear operator $T(g)$ such that,

(i) $T(g_1 g_2) = T(g_1) T(g_2)$,
(ii) $T(e) = I$,

where $e$ is the identity element of $G$ and $I$ is the unit operator on $\mathcal{H}$, then the assignment $g \to T(g)$ is called a *linear representation* of $G$ by operators $T$ on $\mathcal{H}$. Usually, the word *linear* is omitted when referring to linear representations. The dimension of the representation is the dimension of the linear space $\mathcal{H}$.

A representation is called *projective* if (i) above is replaced by

(i′) $T(g_1 g_2) = \epsilon(g_1, g_2)\, T(g_1) T(g_2), \quad |\epsilon(g_1, g_1)| = 1$.

A representation is called *unitary* if each $T(g)$ is a unitary operator.

Two representations $T(g)$ and $Q(g)$ on linear spaces $\mathcal{H}$ and $\mathcal{K}$ are said to be *equivalent* if there exists a linear operator $V$, mapping $\mathcal{H}$ into $\mathcal{K}$ with an inverse $V^{-1}$ such that $Q(g) = V T(g) V^{-1}$.

A subspace $\mathcal{H}_1$ of the space $\mathcal{H}$ of the representation $T(g)$ is called *invariant* if, for $\psi \in \mathcal{H}_1, T(g)\psi \in \mathcal{H}_1$ for all $g \in G$. For every representation there are two *trivial* invariant subspaces, namely the whole space and the null subspace. If a representation $T(g)$ possess only trivial invariant subspaces, it is called *irreducible*.

We shall be concerned with *irreducible, projective, unitary* representations of two particular groups.

#### *1.4.2. Transformation groups*

By a transformation of a set $\Omega$, we mean a one-to-one mapping of the set onto itself. Let $G$ be some group. $G$ is a *transformation group* of the set $\Omega$ if, with each element $g$ of this group, we can associate a transformation $\omega \to g\omega$ in $\Omega$, where for any $\omega \in \Omega$,

(i) $(g_1 g_2)\omega = g_1(g_2 \omega)$ and
(ii) $e\omega = \omega$.

A transformation group $G$ on the set $\Omega$ is called *effective* if the only element $g$ for which $g\omega = \omega$ for all $\omega \in \Omega$, is the identity element $e$ of $G$. An effective group $G$ is called *transitive* on the set $\Omega$, if for any two elements $\omega_1, \omega_2 \in \Omega$, there is some $g \in G$, such that $\omega_2 = g\omega_1$. If $G$ is transitive on a set $\Omega$, then $\Omega$ is called a *homogeneous space* for the group $G$.

For example, the rotation group in three-dimensional Euclidean space is not transitive. A point on a given sphere cannot be changed to a point on a sphere of a different radius by a rotation. However, the unit two-sphere is a homogeneous space for the rotation group.

Let $G$ be a transitive transformation group of a set $\Omega$ and let $\omega_o$ be some fixed point of this set. Let $H$ be the subgroup of elements of $G$ which leave the point $\omega_o$ fixed. $H$ is called the *stationary subgroup* of the point $\omega_o$. Let $\omega_1$ be another point in $\Omega$ and let the transformation $g$ carry $\omega_o$ into $\omega_1$. Then transformations of the form $ghg^{-1}, h \in H$ leave the point $\omega_1$ fixed. The stationary subgroups of the two points are *conjugate* to each other.

Take one of the mutually conjugate stationary subgroups $H$ and denote by $G/H$ the space of left cosets of $G$ with respect to $H$. The set $G/H$ is a homogeneous space for $G$ as a transformation group. There is a one-to-one correspondence between the homogeneous space $\Omega$ and the coset space $G/H$. For example, consider the rotation group in three-dimensional space represented by the group of special orthogonal real $3 \times 3$ matrices $SO(3)$. The set of left cosets $SO(3)/SO(2)$ can be put into one-to-one correspondence with the unit two-sphere.

*1.4.3. Lie algebras and Lie groups*

An abstract Lie algebra $\mathcal{G}$ over the complex or real field is a vector space together with a product $[X, Y]$ such that for all vectors $X, Y, Z$ in $\mathcal{G}$ and $a, b$ in the field,

(i) $[X, Y] = -[Y, X]$,
(ii) $[aX + bY, Z] = a[X, Z] + b[Y, Z]$,
(iii) $[[X, Y], Z] + [[Z, X], Y] + [[Y, Z], X] = 0$.

A representation of an abstract Lie algebra by linear operators on a vector space such as a Hilbert space $\mathcal{H}$, is an algebra homomorphism in the sense that the representation operators have the same product properties as those of the original abstract algebra. For an associative vector space of linear operators, the product operation $[A, B]$ is the commutation operation $[A, B] = AB - BA$.

We will consider representations of two Lie algebras of dimension three. If the basis elements are linear operators $E_1, E_2, E_3$, we may indicate a general element as a linear combination $X = aE_1 + bE_2 + cE_3$. The scalar parameters $\{a, b, c\}$ or a subset of them will then become the parameters of the associated probability distribution.

A Lie group is a topological group which is an analytic manifold. The tangent space to that manifold at the identity of the group is called the Lie algebra of the group. It can be shown that the Lie algebra of a Lie group is an abstract Lie algebra. In the case of a linear representation of a Lie group, the associated Lie algebra can be computed explicitly by differentiating curves through the identity. On the other hand, a (local) Lie group associated with a given Lie algebra can be computed explicitly by the so-called exponential map. (See, for example, Miller (1972).)

For our purposes, we focus upon the parameter space which a Lie group representation inherits from its Lie algebra representation via the exponential map.

*1.5. Families of probability distributions*

Let $G$ be a group and $g \to U(g)$ be an irreducible projective unitary representation of $G$ in the Hilbert space $\mathcal{H}$. For fixed unit vector $\psi_0$, the action of $U(g)$ on $\psi_0$ is designated by,

$$(1.5.1) \qquad \psi_0 \to U(g)\psi_0 = \psi_g.$$

Since each $U(g)$ is unitary, each $\psi_g$ is a unit vector and so can serve as a state. This method of generating a family leads to states designated as "coherent states". The name originated in the field of quantum mechanics. However we use it in a purely group theoretic context, as in Perelomov (1986) and Ali, Antoine, Gazeau (2000). In the two families of probability distributions that we consider in detail, the corresponding families of states are coherent states. See the references given above for properties of coherent states along with examples and suggestions for generalizations. Families of states lead to families of probability distributions. Thus, for self-adjoint operator $O$ in $O$-diagonal Hilbert space $\mathcal{H}$,

$$\mathcal{P}_{\psi_g}\{O \text{ result in } B\} = (U(g)\psi_0, \mathcal{E}(B)U(g)\psi_0) = \sum |(U(g)\psi_0, \eta_i)|^2.$$

In the discrete case, where $\eta_i$ is the eigenvector corresponding the $i$th eigenvalue of $O$, each eigenspace is one-dimensional, and the sum is over all eigenvalues in $B$, and in the continuous case,

$$\mathcal{P}_{\psi_g}\{O \text{ result in } B\} = (U(g)\psi_0, \mathcal{E}(B)U(g)\psi_0) = \int_B |U(g)\psi_0(x)|^2 \mu(dx).$$

## 2. The Poisson family

We construct the Poisson family by first constructing a particular family of coherent states of the form (1.5.1) in an $\ell^2$ Hilbert space $\mathcal{H}_N$. The family is indexed by a parameter set which also indexes a homogeneous space for a certain transformation group; namely, the Weyl–Heisenberg group, denoted $G_W$. Representation operators $T(g), g \in G_W$, acting on a fixed vector in $\mathcal{H}_N$ as in (1.5.1), generate the coherent states which, in turn, generate the family of probability distributions which leads to the Poisson distribution. This provides a context for an inferred probability distribution on the parameter space.

### 2.1. Representations of the Weyl–Heisenberg Lie group

The group $G_W$ can be described abstractly as a three-parameter group with elements $g(s; x_1, x_2)$, for real parameters $s, x_1$ and $x_2$, where the multiplication law is given by

$$(s; x_1, x_2)(t; y_1, y_2) = \left(s + t + \frac{1}{2}(x_1 y_2 - y_1 x_2); \; x_1 + y_1, x_2 + y_2\right).$$

Alternatively, we may consider one real parameter $s$ and one complex parameter $\alpha$, where

(2.1.1) $$\alpha = \frac{1}{\sqrt{2}}(-x_1 + ix_2).$$

Then, $(s; \alpha)(t; \beta) = (s + t + Im(\alpha\beta^*); \; \alpha + \beta)$. The Lie algebra $\mathcal{G}_W$ of the Lie group $G_W$ is a nilpotent three-dimensional algebra. Basis elements can be designated abstractly as $e_1, e_2, e_3$, with commutation relations $[e_1, e_2] = e_3$, $[e_1, e_3] = [e_2, e_3] = 0$. We consider a linear representation of the algebra with basis given by linear operators $E_j$, for $j = 1, 2, 3$, which operate in a Hilbert space $\mathcal{H}$. These operators are such that operators $iE_j$ are self-adjoint with respect to the inner product in $\mathcal{H}$. That property is necessary in order that from this algebra we may construct group representation operators which are unitary. Since the representation is an algebra homomorphism, the linear operators $E_j$ have the same commutation relations as the abstract Lie algebra above.

It will prove to be convenient to consider an alternative basis for the three-dimensional linear space of those operators. Put

$$A = \frac{1}{\sqrt{2}}(E_1 - iE_2), \quad A^\dagger = -\frac{1}{\sqrt{2}}(E_1 + iE_2), \quad I = -iE_3.$$

Note that, due to the fact that the $iE_j$ operators are self-adjoint, $A^\dagger$ is indeed an adjoint operator to $A$. Although $A$ and $A^\dagger$ are not self-adjoint operators, the operator $N = A^\dagger A$ is self-adjoint. We have

(2.1.2) $$[A, A^\dagger] = I, \quad [A, I] = [A^\dagger, I] = 0.$$

A general element of the Lie algebra of operators described above is given by a linear combination of the basis vectors such as $X = isI + \alpha A^\dagger - \alpha^* A$. This form of linear combination derives from $X = x_1 E_1 + x_2 E_2 + sE_3$ and (2.1.1). We may now proceed to obtain a group representation by unitary operators $T$ in the Hilbert space $\mathcal{H}$. By virtue of the exponential map we have, $(s; \alpha) \to T(s; \alpha) = \exp(X)$. Since the $I$ operator commutes with $A$ and $A^\dagger$, we may write $T(s; \alpha) = e^{is} D(\alpha)$ where $D(\alpha) = \exp(\alpha A^\dagger - \alpha^* A)$. It is known that this representation is irreducible.

### 2.2. The Hilbert space of the irreducible representation

The linear operators mentioned above act in a Hilbert space which has been designated abstractly as $\mathcal{H}$. In order to consider concrete formulas for probability distributions, it is necessary to give a concrete form to the Hilbert space. In the case of the Poisson family, the space designated $\mathcal{H}_N$ is realized as an $\ell^2$ space of complex-valued square summable sequences with basis consisting of the eigenvectors of the self-adjoint operator $N$.

By the so-called "ladder" method, using (2.1.2), it has been found that $N$ has a simple discrete spectrum of non-negative integers. Thus, by the general theory, its eigenvectors, $\{\phi_k\}, k = 0, 1, 2, \ldots$, form a complete orthonormal set in $\mathcal{H}$ which forms a basis for the $\ell^2$ Hilbert space realization $\mathcal{H}_N$.

In $\mathcal{H}_N$, we have the following useful properties of the $A$ (annihilation), $A^\dagger$ (creation), and $N$ (number) operators.

(2.2.1)
$$\begin{aligned} A\phi_0 &= 0, \quad A\phi_k = \sqrt{k}\,\phi_{k-1} \quad &&for \quad k = 1, 2, 3, \ldots, \\ A^\dagger \phi_k &= \sqrt{k+1}\,\phi_{k+1} \quad &&for \quad k = 0, 1, 2, 3, \ldots, \\ N\phi_k &= k\phi_k \quad &&for \quad k = 0, 1, 2, 3, \ldots. \end{aligned}$$

Then we can relate $\phi_k$ to $\phi_0$ by

(2.2.2) $$\phi_k = \frac{1}{\sqrt{k!}}(A^\dagger)^k \phi_0 \quad for \quad k = 0, 1, 2, \ldots.$$

### 2.3. Family of coherent states generated by group operators

To construct a family of coherent states in $\mathcal{H}_N$ leading to the Poisson distribution, we operate on the basis vector $\phi_0$ with $D(\alpha)$ operators indexed by complex number $\alpha$, writing

(2.3.1) $$v(\alpha) = D(\alpha)\phi_0.$$

To find an explicit formula for the $v$ vectors, write $D(\alpha)$ as a product of exponential operators. Since $A$ and $A^\dagger$ do not commute, we do not have the property which pertains to scalar exponential functions that $D(\alpha) = \exp(\alpha A^\dagger)\exp(-\alpha^* A)$. We use the Baker-Campbell-Hausdorff operator identity:

$$(2.3.2) \qquad \exp(O_1)\exp(O_2) = \exp\left(\frac{1}{2}[O_1, O_2]\right)\exp(O_1 + O_2)$$

which is valid in the case where the commutator $[O_1, O_2]$ commutes with both operators $O_1$ and $O_2$. Putting $O_1 = \alpha A^\dagger, O_2 = -\alpha^* A$, we have

$$D(\alpha)\phi_0 = e^{-|\alpha|^2/2}\exp(\alpha A^\dagger)\exp(-\alpha^* A)\ \phi_0.$$

For linear operators, we have the same kind of expansion for an exponential operator function as for a scalar function. Expanding $\exp(-\alpha^* A)\phi_0$ and using (2.2.1), we find that $\exp(-\alpha^* A)\phi_0 = I\ \phi_0$. Then from (2.2.2), we see that $(A^\dagger)^k \phi_0 = \sqrt{k!}\ \phi_k$. From (2.2.1) and (2.3.1),

$$(2.3.3) \qquad D(\alpha)\phi_0 = e^{-|\alpha|^2/2}\sum_{k=0}^{\infty}\frac{\alpha^k}{\sqrt{k!}}\ \phi_k.$$

### 2.4. Family of probability distributions

Let the observable (self-adjoint) number operator $N$ represent the physical quantity being "counted" with possible outcomes $0, 1, 2, \ldots$. Using the family of coherent states given above, the probability distributions are,

$$\mathcal{P}_{v(\alpha)}\{\text{result} = n\} = |(\phi_n, v(\alpha))|^2.$$

By expression (2.3.3) and the orthogonality of the basis vectors,

$$(2.4.1) \qquad (\phi_n, v(\alpha)) = e^{-|\alpha|^2/2}\sum_{k=0}^{\infty}\frac{\alpha^k}{\sqrt{k!}}(\phi_n, \phi_k) = e^{-|\alpha|^2}\frac{\left(|\alpha|^2\right)^n}{n!}$$

for $n = 0, 1, 2, 3, \ldots$. Taking the modulus squared, using (1.2.1), we have the formula for the Poisson family,

$$\mathcal{P}_{v(\alpha)}\{\text{result} = n\} = e^{-|\alpha|^2}\frac{\left(|\alpha|^2\right)^n}{n!} \quad for \quad n = 0, 1, 2, 3, \ldots.$$

Put the Poisson parameter $\lambda = |\alpha|^2$. Thus we see that $\lambda$ is real and nonnegative.

It may be remarked that this is a complicated method for obtaining the Poisson family. The point is that we now have a context in which to infer a probability distribution on the parameter space, given an observed Poisson value $n$.

### 2.5. POV measures versus PV measures

Consider the definition of a projection-valued (PV) measure, or *spectral measure* (see, for example, Busch, Grabowski, and Lahti (1995)), which had been introduced heuristically in Section 1.2.

**Definition.** Let $\mathcal{B}(\mathbb{R})$ denote the Borel sets of the real line $\mathbb{R}$ and $\Lambda(\mathcal{H})$ denote the set of bounded linear operators on $\mathcal{H}$. A mapping $\mathcal{E}: \mathcal{B}(\mathbb{R}) \to \Lambda(\mathcal{H})$ is a projection valued (PV) measure, or a spectral measure, if

$$(2.5.1) \begin{aligned} &\mathcal{E}(B) = \mathcal{E}^\dagger(B) = \mathcal{E}^2(B), \qquad \forall B \in \mathcal{B}(\mathbb{R}), \\ &\mathcal{E}(\mathbb{R}) = I, \;\; \mathcal{E}(\cup_i B_i) = \sum_i \mathcal{E}(B_i) \;\; \text{for all disjoint sequences } \{B_i\} \subset \mathcal{B}(\mathbb{R}), \end{aligned}$$

where the series converges in the weak operator topology. This spectral measure gives rise to the definition of a unique self-adjoint operator $O$ defined by $O = \int x \mathcal{E}(dx)$, with its domain of definition

$$(2.5.2) \;\; \mathcal{D}(O) = \left\{ \psi \in \mathcal{H}, s.t. \left(\psi, \int x^2 \mathcal{E}(dx)\psi\right) = \int x^2 (\psi, \mathcal{E}(dx)\psi) \;\; converges \right\}.$$

Given a self-adjoint operator $O$ with domain of definition $\mathcal{D}(O) \subseteq \mathcal{H}$, there is a unique PV measure $\mathcal{E}: \mathcal{B}(\mathbb{R}) \to \Lambda(\mathcal{H})$ such that $\mathcal{D}(O)$ is (2.5.2), and for any $\psi \in \mathcal{D}(O)$, $(\psi, O\psi) = \int x (\psi, \mathcal{E}(dx)\psi)$. Therefore there is a one-to-one correspondence between self-adjoint operators $O$ and real PV measures $\mathcal{E}$.

In the case of the Poisson distribution, the self-adjoint operator is $N$ with normalized eigenvectors $\{\phi_k\}$ as the orthonormal basis of $\mathcal{H}$, and $sp(N) = \{0, 1, 2, \ldots\}$. For an inferential probability measure operator on the parameter space associated with the Poisson distribution, we will have neither a PV measure nor a self-adjoint operator. Instead we will have a *subnormal* operator and a so-called positive operator-valued (POV) measure, where the first line of (2.5.1) is amended to read

$$\mathcal{E}(B) \text{ is a positive operator for all } B \in \mathcal{B}(\mathbb{R}).$$

The operators for PV measures are *projections*. The properties prescribed for those projections $\mathcal{E}(B)$ are just those needed so that the corresponding inner products $(\psi, \mathcal{E}(dx)\psi)$ for vector states $\psi$ will have the properties of probabilities. However, for the definition of probability there is no requirement that the operators be projections. In fact, if they are *positive* operators, they can still lead to probabilities when put together with states.

### 2.6. An invariant measure on the parameter space

Now we consider the inferential case. In a sense we reverse the roles of states and observable operators. If the Poisson value $n$ was observed, what was formerly a vector $\phi_n$ denoting a one-dimensional projection $\mathcal{E}(\{n\})$, now becomes a *state*. What was formerly a family of coherent sates, $v(\alpha)$, now leads to the construction of a POV measure on the parameter set.

In order to obtain a measure on the parameter space $\mathbb{C}$ which is invariant to an operator $D(\beta)$, for arbitrary complex number $\beta$, we need to see how the operator transforms a coherent state $v(\alpha)$. Consider $D(\beta)v(\alpha) = D(\beta)D(\alpha)\phi_0$, Using the Baker–Campbell–Hausdorff identity (2.3.2) with $O_1 = \beta A^\dagger - \beta^* A$ and $O_2 = \alpha A^\dagger - \alpha^* A$, we have $D(\beta)D(\alpha)\phi_0 = e^{iIm(\beta\alpha^*)}D(\beta+\alpha)\phi_0$, As *states*, $D(\beta)v(\alpha) = v(\beta+\alpha)$. Thus, the operator $D$ acts as a translation operator on the complex plane so that the invariant measure, $d\mu(\alpha)$, is just Lebesgue measure, $d\mu(\alpha) = c d\alpha_1 d\alpha_2$, where $\alpha = \alpha_1 + i\alpha_2$, and where $c = 1/\pi$ by normalization.

## 2.7. An inferred distribution on the parameter space

By general group theory, the irreducibility of the group representation by unitary operators in the Hilbert space $\mathcal{H}$ implies that the coherent states are complete in $\mathcal{H}$. (See, for example, Perelomov (1986)). Thus, for any vectors $\psi_1$ and $\psi_2$ in $\mathcal{H}$,

$$(\psi_1, \psi_2) = \int (\psi_1, v(\alpha))(v(\alpha), \psi_2) \, d\mu(\alpha).$$

The coherent states form a so-called "overcomplete" basis for a Hilbert space in the sense that they are complete and can be normalized but are not orthogonal. The Hilbert space $\mathcal{H}_{CS}$ which they span may be visualized as a proper subspace of $L^2(\mathbb{C})$, the space of square integrable functions of a complex variable with inner product $(f(\alpha), g(\alpha))_{CS} = \int_{\mathbb{C}} f(\alpha)^* g(\alpha) d\mu(\alpha)$. In Ali, Antoine, Gazeau(2000), an isometric map $\rho$ is given which associates an element $\phi$ in $\mathcal{H}_N$ with element (function of $\alpha$) in $\mathcal{H}_{CS}$: $\rho(\phi) = (\phi, v(\alpha))_{H_N}$. In $\mathcal{H}_{CS}$ we construct a POV measure $M$ which leads to a probability distribution for $\alpha$ defined by

$$\mathcal{P}\{\alpha \in \Delta \ for \ state \ \psi\} = (\psi, M(\Delta)\psi) = \int_\Delta |(\psi, v(\alpha))|^2 d\mu(\alpha)$$

for complex Borel set $\Delta$ and for state $\psi$. In particular, consider the (eigenvector) state $\psi = \phi_n$ corresponding to an observed Poisson value $n$, an eigenvalue of the self-adjoint number operator $N$. $\mathcal{P}\{\alpha \in \Delta \ for \ state \ \phi_n\} = \int_\Delta |(\phi_n, v(\alpha))|^2 d\mu(\alpha)$. This provides us with a probability distribution on the whole parameter space, namely, the complex plane. But the Poisson parameter, a real number, is the modulus squared $|\alpha|^2$ of $\alpha$. Expressing $\alpha$ in polar coordinates, $\alpha = re^{i\theta}$, with $r > 0$ and $0 \leq \theta < 2\pi$, we obtain the invariant measure $d\mu(\alpha) = \frac{1}{2\pi} dr^2 d\theta = \frac{1}{\pi} r dr d\theta$. Then integrating $\theta$ from 0 to $2\pi$, we obtain the marginal distribution for $r^2$ as follows. For real Borel set $B$,

$$\mathcal{P}\{r^2 \in B \ for \ state \ \phi_n\} = \int_B e^{-r^2} \frac{(r^2)^n}{n!} \left(\frac{1}{2\pi} \int_0^{2\pi} d\theta\right) dr^2 = \int_B e^{-\lambda} \frac{\lambda^n}{n!} \, d\lambda,$$

where $\lambda = r^2$ and the expression for $|(\phi_n, w(\alpha))|^2$ is obtained similarly as in (2.4.1). We see that this corresponds to a Bayes posterior distribution with uniform prior distribution for the parameter $\lambda$.

## 3. The binomial family

We construct the binomial family similarly as was done for the Poisson family. In this case the coherent states are built from irreducible representations of the Lie algebra of the rotation group $SO(3)$ of real $3 \times 3$ orthogonal matrices with determinant one, instead of the Weyl–Heisenberg group. The Weyl–Heisenberg Lie algebra is three-dimensional nilpotent whereas the Lie algebra corresponding to $SO(3)$ is three-dimensional simple.

## 3.1. The rotation group and Lie algebra

Although there are nine matrix elements in a $3 \times 3$ real matrix, the constraints of orthogonality and unit determinant for an element $g$ of $SO(3)$, imply that $g$ can be identified by three real parameters. There are two general methods for indicating the parameters. One way is by using the three Euler angles. The other way is to specify an axis and an angle of rotation about that axis.

The rotation group is locally isomorphic to the group $SU(2)$ of $2 \times 2$ complex unitary matrices of determinant one. An element $u$ of $SU(2)$ is identified by two complex numbers $\alpha$ and $\beta$ where $|\alpha|^2 + |\beta|^2 = 1$. The relationship between $SO(3)$ and $SU(2)$ is that of a unit sphere to its stereographic projection upon the complex plane as shown in Naimark (1964). Although the relationship is actually homomorphic, (one $g$ to two $u$), they have the same Lie algebra and so can be used interchangeably in the context presented here. It is more intuitive to work with $SO(3)$ but from the point of view of the binomial distribution, it will turn out to be more pertinent to work with $SU(2)$. Both $SO(3)$ and $SU(2)$ are compact as topological groups (Vilenkin (1968)).

In this case, since we start with matrices, basis elements of the Lie algebra can be easily obtained by differentiating the three matrices corresponding to the subgroups of rotations about the three spatial coordinate axes $(x, y, z)$. Thus, for example, the subgroup of $SO(3)$, indicating rotations about the $z$ axis is given by

$$a_3(t) = \begin{pmatrix} \cos t & -\sin t & 0 \\ \sin t & \cos t & 0 \\ 0 & 0 & 1 \end{pmatrix}.$$

Differentiating each matrix element with respect to $t$ and then setting $t = 0$, we obtain the algebra basis element

$$e_3 = \begin{pmatrix} 0 & -1 & 0 \\ 1 & 0 & 0 \\ 0 & 0 & 0 \end{pmatrix}.$$

Similarly, we obtain the three basis elements, $e_1, e_2, e_3$ with commutation relations

(3.1.1) $\qquad [e_1, e_2] = e_3, \ [e_2, e_3] = e_1, \ [e_3, e_1] = e_2.$

## 3.2. A homogeneous space for the group

The rotation group $G = SO(3)$ acts as a transformation group in three-dimensional Euclidean space. However, $SO(3)$ is not transitive on the whole space. It is transitive on spheres. We take the unit two-sphere $\mathbb{S}^2$ as a homogeneous space for the group. But there is not a one-to-one relationship between group elements and points on the unit sphere. A one-to-one relationship (excluding the South Pole of the sphere) is provided by the cosets $G/H_{NP}$ of the group with respect to the isotropy subgroup $H_{NP}$ of the North Pole $(0, 0, 1)$ of $\mathbb{S}^2$. In $SO(3)$, the subgroup is the group $a_3(t)$ of rotations about the $z$ axis. In $SU(2)$, the subgroup $U(1)$ is the set of diagonal matrices $h(t)$ with diagonal elements $e^{it}$ and $e^{-it}$.

Following Perelomov (1986), we consider cosets $SO(3)/SO(2)$ or cosets $SU(2)/U(1)$. The one-to-one relationship of cosets $SU(2)/U(1)$ with the unit sphere $\mathbb{S}^2$ (excluding the South Pole) is given in the following manner. Given a point $\nu$ on the unit sphere indicated by

(3.2.1) $\qquad\qquad \nu = (\sin\theta\cos\gamma, \ \sin\theta\sin\gamma, \ \cos\theta),$

associate the coset $g_\nu$ where, $g_\nu = \exp\left(\frac{i\theta}{2}(\sin\gamma M_1 - \cos\gamma M_2)\right)$, where the matrices $M_1$ and $M_2$ are given by $M_1 = \begin{pmatrix} 0 & 1 \\ 1 & 0 \end{pmatrix}$, $M_2 = \begin{pmatrix} 0 & -i \\ i & 0 \end{pmatrix}$. In terms of rotations, the matrix describes a rotation by angle $\theta$ about the axis indicated by direction $(\sin\gamma, -\cos\gamma, 0)$ which is perpendicular to both the North Pole and the vector $\nu$. We can express a general element $u$ of $SU(2)$ by

(3.2.2)  $\quad u = g_\nu h, \quad where \ g_\nu \in \ coset \ SU(2)/U(1), \quad h \in U(1).$

### 3.3. Irreducible representations

Now we consider a representation of the algebra with basis given by linear operators $E_k$, for $k = 1, 2, 3$ which operate in a Hilbert space $\mathcal{H}$. Since the group is compact, general theory provides the result that irreducible representations correspond to finite dimensional Hilbert spaces. In the algebra representation space, the basis elements have the same commutation relations as (3.1.1). Also, we require that the operators $J_k = iE_k$ be self-adjoint with respect to the inner product of $\mathcal{H}$.

Similarly as in Section 2.1, introduce creation and annihilation operators $J_+ = J_1 + iJ_2$, $J_- = J_1 - iJ_2$. It is known that the complete set of irreducible representations of the Lie algebra is indexed by a non-negative integer or half-integer $j$ while the dimension of the representation space is $2j+1$. Correspondingly, the complete set of unitary irreducible representations $T_j(u)$ of the group $SU(2)$ is given by $j = 0, 1/2, 1, 3/2, \ldots$. Since the relationship of $SU(2)$ to the rotation group $SO(3)$ is two-to-one, the irreducible representations of the rotations group are more properly indexed by the non-negative integers, omitting the half-integers. We will see that the parameter $n$ for the binomial distribution is equal to $2j$ implying that we need the non-negative half-integers in the list. Thus we focus upon the group $SU(2)$. Choose and fix the number $j$ indicating a definite Hilbert space $\mathcal{H}_j$. An orthonormal basis for $\mathcal{H}_j$ is provided by the eigenvectors $\phi_m$ of the self-adjoint operator $J_3$ which, for fixed $j$, has a simple discrete and finite spectrum indexed by $m = -j, -j+1, \ldots, j-1, j$. As operators in $\mathcal{H}_j$, $J_+, J_-$ and $J_3$ have creation, annihilation, and number properties similarly as in (2.2.1):

$$J_+\phi_j = 0, \quad J_+\phi_m = \sqrt{(j-m)(j+m+1)} \ \phi_{m+1}, \quad for \ m = -j, -j+1, \ldots, j-1,$$
(3.3.1) $J_-\phi_{-j} = 0, \quad J_-\phi_m = \sqrt{(j+m)(j-m+1)} \ \phi_{m-1}, \quad for \ m = -j+1, -j+2, \ldots, j,$
$$J_3\phi_m = m\phi_m, \qquad\qquad for \ m = -j, -j+1, \ldots, j-1, j.$$

Note that in (2.2.1) there is a minimum basis vector $\phi_0$, but no maximum indicating an infinite dimensional Hilbert space. Here we have both a minimum and a maximum basis vector. We relate $\phi_m$ and $\phi_{-j}$ by

(3.3.2) $$\phi_m = \sqrt{\frac{(j-m)!}{(j+m)!(2j)!}} \ (J_+)^{j+m} \phi_{-j}.$$

For fixed number $j$, and for each $u \in SU(2)$, let $u \to T^j(u)$ denote an irreducible representation of $SU(2)$ in the Hilbert space $\mathcal{H}_j$, where each operator $T^j$ is unitary with respect to the inner product in $\mathcal{H}_j$. From (3.2.2), we have $T^j(u) = D(\nu)T^j(h)$ for $h \in U(1)$, where $D(\nu) = T^j(g_\nu)$,

$$T^j(g_\nu) = \exp(i\theta(\sin\gamma J_1 - \cos\gamma J_2)), \quad for \ 0 \leq \theta < \pi.$$

It can be shown that for $h \in U(1)$, $T^j(h)$ is a scalar factor of modulus one. Thus we focus upon $D(\nu)$.

### 3.4. The family of coherent states and the binomial family

We choose the element $\phi_{-j}$ as the fixed vector in $\mathcal{H}_j$. Then similarly as in Section 2, the family of coherent states, is given by

$$w(\nu) = D(\nu)\phi_{-j} = \exp\left(i\theta(\sin\gamma J_1 - \cos\gamma J_2)\right)\phi_{-j}, \quad for\ 0 \le \theta < \pi.$$

As in Perelomov(1986), we find it convenient to re-parameterize, similarly as in Section 2.1, and use one complex parameter $\xi$ along with the creation and annihilation operators instead of the two real angle parameters with the $J_1$ and $J_2$ operators. Thus for $\xi = \frac{i\theta}{2}(\sin\gamma + i\cos\gamma)$, we have $D(\xi) = \exp(\xi J_+ - \xi^* J_-)$. We seek an explicit expression for the coherent states $w(\nu)$. As in Section 2.3, the method is to factor the exponential function. The Baker–Campbell Hausdorff formula is not convenient to use in this case as it was in Section 2.3. Instead, the Gauss decomposition of the group $SL(2,\mathbb{C})$ is used. We obtain, $D(\xi) = \exp(\zeta J_+)\exp(\eta J_3)\exp(\zeta' J_-)$, where $\zeta = -\tan(\theta/2)e^{-i\gamma}$, $\eta = -2\ln\cos|\xi| = \ln(1+|\zeta|^2)$, $\zeta' = -\zeta^*$. Finally, using (3.3.1) and (3.3.2), we obtain coherent states

$$w(\zeta) = \sum_{m=-j}^{j} \sqrt{\frac{(2j)!}{(j+m)!(j-m)!}} \frac{\zeta^{j+m}}{(1+|\zeta|^2)^j} \phi_m.$$

In terms of angle parameters,

$$(3.4.1)\quad w(\theta,\gamma) = \sum_{m=-j}^{j} \sqrt{\frac{(2j)!}{(j+m)!(j-m)!}} \left(-\sin\frac{\theta}{2}\right)^{j+m} \left(\cos\frac{\theta}{2}\right)^{j-m} e^{-i(j+m)\gamma} \phi_m.$$

Thus, noting that the possible result values are the eigenvalues of $J_3$, namely, $m = -j, -j+1, \ldots, j-1, j$, we have

$$\mathcal{P}\{result = \ell,\ \text{when the state is } w(\theta,\gamma)\} = |(\phi_\ell, w(\theta,\gamma))|^2.$$

Using the fact that the eigenvectors $\phi_m$ of $J_3$ are orthonormal, we have,

$$(3.4.2)\quad (\phi_\ell, w(\theta,\gamma)) = \sqrt{\frac{(2j)!}{(j+\ell)!(j-\ell)!}} \left(-\sin\frac{\theta}{2}\right)^{j+\ell} \left(\cos\frac{\theta}{2}\right)^{j-\ell} e^{-i(j+\ell)\gamma}.$$

Upon taking the modulus squared, we have,

$$\mathcal{P}_{w(\theta,\gamma)}\{result = \ell\} = \frac{(2j)!}{(j+\ell)!(j-\ell)!}\left(\sin^2\frac{\theta}{2}\right)^{j+\ell}\left(\cos^2\frac{\theta}{2}\right)^{j-\ell},$$

for $\ell = -j, -j+1, \ldots, j-1, j$. For the binomial $(n,p)$ distribution, put $n = 2j$, renumber the possible values by $k = j + l$ and put $p = \sin^2\theta/2$.

### 3.5. An inferred distribution on the parameter space

The parameters $\theta$ and $\gamma$ index the parameter space, that is, the points of the unit two-sphere, which is isomorphic to the cosets $SU(2)/U(1)$, or equivalently $SO(3)/SO(2)$, and where points are given by the (three-dimensional) vector $\nu$ as in (3.2.1). In other words, we can take the unit sphere to be the parameter space. The coherent states, also indexed by the point $\nu$, are complete in the Hilbert space $\mathcal{H}_j$. As before, we have an isometric map from $\mathcal{H}_j$ to the Hilbert space $\mathcal{H}_{CS}^j$ spanned by the coherent states. Since $D(\nu)$ takes one coherent state into another coherent state, we have the action of $D(\nu)$ on $\mathcal{H}_{CS}^j$. The (normalized) measure invariant to the action of $D(\nu)$ is Lebesgue measure on the sphere: $d(\theta,\gamma) = \frac{2j+1}{4\pi}\sin\theta d\theta d\gamma$.

Suppose that we have an observed binomial count value $\tilde{k}$ which, with $\tilde{k} = j + \tilde{\ell}$, gives $\tilde{\ell} = \tilde{k} - j$, for possible values $\tilde{\ell} = -j, -j+1, \ldots, j-1, j$ corresponding to possible values $\tilde{k} = 0, 1, 2, \ldots, 2j$. Then the corresponding inferred distribution on the parameter space, derived from a POV measure is

$$\mathcal{P}\{(\theta, \gamma) \in \Delta \text{ when the state is } \phi_{\tilde{\ell}}\} = \int_\Delta |(\phi_{\tilde{\ell}}, w(\nu))|^2 d(\theta, \gamma),$$

where the inner product inside the integral sign is that of $\mathcal{H}_j$. From the expression for $w(\nu)$ in (3.4.1) and the fact that the vectors $\phi_m$ are orthonormal in $\mathcal{H}_j$ giving the inner product (3.4.2), we have the joint distribution of $\theta$ and $\gamma$:

$$\frac{2j+1}{4\pi} \iint_\Delta \left| \sqrt{\frac{(2j)!}{(j+\tilde{\ell})!(j-\tilde{\ell})!}} \left(-\sin\frac{\theta}{2}\right)^{j+\tilde{\ell}} \left(\cos\frac{\theta}{2}\right)^{j-\tilde{\ell}} e^{-i(j+\tilde{\ell})\gamma} \right|^2 \sin\theta d\theta d\gamma,$$

$$= \frac{2j+1}{4\pi} \iint_\Delta \frac{(2j)!}{(j+\tilde{\ell})!(j-\tilde{\ell})!} \left(\sin^2\frac{\theta}{2}\right)^{j+\tilde{\ell}} \left(\cos^2\frac{\theta}{2}\right)^{j-\tilde{\ell}} \sin\theta d\theta d\gamma.$$

For the marginal distribution of $\theta$, we integrate $\gamma$ from $0 \leq \gamma < 2\pi$ obtaining. For $B$, a Borel set in $[0, \pi)$,

$$\mathcal{P}\{\theta \in B \text{ when the state is } \phi_{\tilde{\ell}}\}$$

$$= \frac{2j+1}{2} \int_B \frac{(2j)!}{(j+\tilde{\ell})!(j-\tilde{\ell})!} \left(\sin^2\frac{\theta}{2}\right)^{j+\tilde{\ell}} \left(\cos^2\frac{\theta}{2}\right)^{j-\tilde{\ell}} \sin\theta d\theta,$$

$$= \frac{n+1}{2} \int_B \frac{n!}{\tilde{k}!(n-\tilde{k})!} \left(\sin^2\frac{\theta}{2}\right)^{\tilde{k}} \left(\cos^2\frac{\theta}{2}\right)^{n-\tilde{k}} \sin\theta d\theta.$$

$$= (n+1) \int_B \frac{n!}{\tilde{k}!(n-\tilde{k})!} p^{\tilde{k}} (1-p)^{n-\tilde{k}} dp,$$

for $p = \sin^2\frac{\theta}{2}$, implying a uniform prior distribution for the canonical parameter $p$.

## 4. Discussion

We have constructed a group theoretic context for the two discrete probability distributions, Poisson and binomial. Similarly as other group invariance methods, the idea is to construct probability families by group action. However, in contrast to others, we have neither a pivotal function nor group action on the value space of the random variable. Thus our method is applicable to the discrete case. In this paper the Poisson and binomial families were constructed by using the properties of certain families of vectors (coherent states) which due to their completeness property enable the construction of measures leading to inferred distributions on the parameter spaces. The formulas for the inferred distributions obtained in those two examples coincided with Bayesian posterior distributions in the case where the prior distributions were uniform. We emphasize the fact that although the formulas for the two methods coincide in the end result, the two methods are distinctly different.

This difference may be illustrated by considering Thomas Bayes' justification for a uniform prior in the binomial case elucidated in Stigler (1982). Here the emphasis is not on the parameter itself, but rather on the *marginal* distribution of the binomial random variable $X$ obtained from the *joint distribution* of $X$ and parameter $p$. Starting with a joint distribution applied to a particular binomial physical situation (billiard table example) in which the parameter has uniform distribution and integrating out the parameter obtaining the marginal distribution for $X$, one obtains the result that $X$ has a discrete uniform distribution. Then the reasoning is that, in the face of no prior knowledge, one assumes a discrete uniform distribution for $X$ for all $n$ implying a uniform prior distribution for $p$. Stigler notes that if $P(X = k)$ is constant, then so is $P(f(X) = f(k))$ for any strictly monotone function $f(x)$, thus answering the objection raised against the principle of insufficient reason where a uniform distribution for a given parameter would not be uniform for every monotone function of it. The argument raised against this approach of Bayes is that it is too restrictive in that it "is very strongly tied to the binomial model."

The group theoretic method operates in a different context. There is no joint distribution of random variable and parameter and consequently no marginal distribution for the random variable. One starts by constructing an ordinary family of probability distributions indexed by parameters obtained from a chosen parametric group. To obtain the inferred distribution on the parameter space, the roles are exchanged in that the observed value of the original random variable acts as a parameter and the former parameters are treated as random. The original random variable and parameters are never random at the same time. The reversal in roles is possible technically because of the completeness property of the coherent states which were used in the first place to construct the family. In the binomial case, the relevant group is the matrix group $SU(2)$ and the consequent invariant distribution is Lebesgue measure on the 2-sphere. Upon integrating out the azimuthal angle $\gamma$, we obtain $sin\theta d\theta$ for the polar angle which, with a slight change of variable yields a probability distribution for parameter $p$ which is the same formula as a Bayesian posterior based upon uniform prior for $p$. We obtained similar results for the Poisson distribution using the Weyl-Heisenberg group. Clearly, the list of discrete distributions with associated groups can be extended. Results (unpublished) relating to matrix group $SU(1, 1)$ have been obtained for the negative binomial distribution. Unlike the case of binomial, our results do not imply uniform prior for the commonly used parameter $p$ as given, for example, in [5].

Efron (1998) has indicated a relationship between the fiducial method of inference and the Bayesian method as follows: "By 'objective Bayes' I mean a Bayesian theory in which the subjective element is removed from the choice of prior distribution; in practical terms a universal recipe for applying Bayes theorem in the absence of prior information. A widely accepted objective Bayes theory, which fiducial inference was intended to be, would be of immense theoretical and practical importance."

From the Bayesian point of view, one may interpret this paper as an objective method for obtaining a reference prior in the absence of prior information. From another point of view, one might interpret this paper as a way of obtaining inferred distributions on parameter spaces without the use of the Bayes method.

**Acknowledgments.** We thank the editor, the associate editor and the referee for their valuable comments and suggestions.